\title{ ~~\\ Counting carefree couples}
\author{Pieter Moree}
\documentclass[12pt]{article}
\usepackage{amssymb, latexsym, amsfonts}
\def\@ptsize{2}
\newtheorem{Thm}{Theorem}
\newtheorem{lem}{Lemma}
\newtheorem{Problem}{Problem}

\newcommand{\qed}{\hfill $\Box$}

\begin{document}
\date{}
\maketitle
\begin{abstract}
\noindent A pair of natural numbers $(a,b)$ such that
$a$ is both squarefree and coprime to $b$ is called a carefree couple. 
A result conjectured by Manfred Schroeder (in his book `Number theory in
science and communication') 
on carefree couples  and a
variant of it are established using standard arguments from elementary analytic
number theory. Also a related conjecture of
Schroeder on triples of integers that are pairwise coprime is proved.

\end{abstract}
\section{Introduction}
It is well known that the probability that an integer is squarefree is
$6/\pi^2$. 
Also the probability that two given integers are coprime
is $6/\pi^2$. 
(More generally the probability that $n$ positive integers chosen 
arbitrarily and independently are coprime is well-known \cite{Lehmer, ny, sylv} to be
$1/\zeta(n)$, where $\zeta$ is Riemann's zeta function. For some generalizations see
e.g. \cite{CS, Coh0, HSM, ny2, Por}.) 
One can wonder how `statistically independent' squarefreeness
and coprimality are.
To this end one could for example consider the probability
that of two random natural
numbers $a$ and $b$, $a$ is both squarefree and
coprime to $b$. Let us call such a couple $(a,b)$ {\it carefree}. If
$b$ is also squarefree, we say that $(a,b)$ is a 
{\it strongly carefree} couple.
Let us denote by $C_1(x)$ the number of carefree couples $(a,b)$ with
both $a\le x$ and $b\le x$ and, similarly, let $C_2(x)$ denote the number
of strongly carefree couples $(a,b)$ with both $a\le x$ and $b\le x$.\\
\indent The purpose of this note is to establish the following result, part
of which was conjectured, on the basis
of heuristic arguments, by Manfred Schroeder \cite[p. 54]{schroed}. (In it and in the
rest of the paper the mathematical symbol $p$ is exclusively used to denote primes.)
\begin{Thm}
\label{main}
We have
\begin{equation}
\label{een}
C_1(x)={x^2\over \zeta(2)}\prod_p\Big(1-{1\over p(p+1)}\Big)
+O(x\log x),
\end{equation}
and
\begin{equation}
\label{twee}
C_2(x)={x^2\over \zeta(2)^2}\prod_p\Big(1-{1\over (p+1)^2}\Big)+O(x^{3/2}).
\end{equation}
\end{Thm}
The interpretation of Theorem \ref{main} is that the probability for
a couple to be carefree is
\begin{equation}
\label{probeen}
K_1:={1\over \zeta(2)}\prod_p\Big(1-{1\over p(p+1)}\Big)\approx 0.42824950567709444022
\end{equation}
and to be strongly carefree is
\begin{equation}
\label{probtwee}
K_2:={1\over \zeta(2)^2}\prod_p\Big(1-{1\over (p+1)^2}\Big)\approx
0.28674742843447873411
\end{equation}
Using the identity $\zeta(n)=\prod_p (1-p^{-n})^{-1}$ valid for $n>1$ we 
can alternatively write
\begin{equation}
\label{trafo}
K_2={1\over \zeta(2)}
\prod_p\Big(1-{2\over p(p+1)}\Big)=
\prod_p\Big(1-{1\over p}\Big)^2\Big(1+{2\over p}\Big).
\end{equation}
For $m\ge 3$ and $0\le k\le m$ we put
\begin{equation}
\label{zkm}
Z_k(m) = \prod_p \Big(1 + {k-1\over p^m} - {k\over p^{m-1}}\Big).
\end{equation}
Note that $Z_2(3)=K_1$ and $Z_3(3)=K_2$.\\
\indent The constants $K_1$ and $K_2$ we could call the {\it carefree}, respectively 
{\it strongly carefree constant}, cf. \cite[Section 2.5]{Finch}.\\
\indent Assuming independence of squarefreeness and coprimality we would expect
that $K_1=\zeta(2)^{-2}$ and $K_2=\zeta(2)^{-3}$.
Now note that 
$$K_1={1\over \zeta(2)^2}\prod_p\Big(1+{1\over (p+1)(p^2-1)}\Big),~K_2={1\over \zeta(2)^3}\prod_p\Big(1+{2p+1\over (p+1)^2(p^2-1)}\Big).$$
We have
$\zeta(2)^2K_1\approx 1.15876$ and $\zeta(2)^3K_2\approx 1.27627$.
Thus, there is a positive correlation between squarefreeness and
coprimality.\\
\indent Let
$I_3(x)$ denote the number of triples $(a,b,c)$ with
$a\le x,~b\le x,c\le x$ such that $(a,b)=(a,c)=(b,c)=1$.
Schroeder \cite[Section 4.4]{schroed}
claims that $I_3(x)\sim K_2x^3$. Indeed, in Section \ref{st2} we will prove the following result.
\begin{Thm}
\label{main2}
We have $I_3(x)=K_2x^3+O(x^2\log^2 x).$
\end{Thm}

\indent The work described in this note was carried out in 2000 and with some improvement in the error
terms was posted on the arXiv in September of 2005 \cite{Morun}, with the remark
that it was not intended for publication in a research journal as the methods used
involve only rather elementary and standard analytic number theory. Over the years 
various authors referred to \cite{Morun}, and this induced me to try to publish it in
a mathematical newsletter. 
(For publications in this area after 2005 see, e.g, \cite{AH,FF1,FF2,FF3,FF4,Hey,Hu-1,Hu,toth1,toth2}.)
In 
\cite{Morun} there was a mistake in the proof of (2) leading to an 
error term of $O(x\log^3 x)$, rather than $O(x^{3/2})$. 
Except for this, the present version has essentially the same mathematical 
content as the earlier one, but is written in a less carefree way and with the
mathematical details more spelled out.

\section{Proofs}
As usual we let $\mu$ denote the M\"obius function and $\varphi$ Euler's totient function. Note
that $n$ is squarefree if and only if $\mu(n)^2=1$.
We will repeatedly make use of the basic identities
\begin{equation}
\label{startmu}
\sum_{d|n}\mu(d)=\cases{1 & if $n=1$;\cr 0 & otherwise,}
\end{equation}
and
\begin{equation}
\label{startphi}
{\varphi(n)\over n}=\sum_{d|n}{\mu(d)\over d}=\prod_{p|n}(1-{1\over p}).
\end{equation}
We will also use several times that if $s$ is a complex number and $f$ a multiplicative function
such that
$\sum_p \sum_{\nu\ge 1}|f(p^{\nu})p^{-\nu s}|<\infty$, then
\begin{equation}
\label{pprod}
\sum_{n=1}^{\infty}{f(n)\over n^s}=\sum_p \sum_{\nu\ge 1}{f(p^{\nu})\over p^{\nu s}}.
\end{equation}
(For a proof see, e.g., Tenenbaum \cite[p. 107]{tenenbaum}.)\\
\indent In the proof of Theorem
\ref{main} we will make use of the following lemma.
\begin{lem}
\label{mainlemma}
Let $d\ge 1$ be arbitrary.
Put
$$S_d(x)=\sum_{n\le x\atop (d,n)=1}\mu(n)^2.$$
We have
\begin{equation}
\label{esdee}
S_d(x)
={x\over \zeta(2)\prod_{p|d}
(1+{1\over p})}+O(2^{\omega(d)}\sqrt{x}),
\end{equation}
where $\omega(d)$ denotes the number of distinct prime divisors of $d$.
\end{lem}
{\it Proof}. Let $T_d(x)$ denote the number of natural numbers
$n\le x$ that are coprime to $d$. Using (\ref{startmu}) and (\ref{startphi}) and 
$[x]=x+O(1)$ we
deduce that
\begin{equation}
\label{teedee}
T_d(x)=\sum_{n\le x\atop (n,d)=1}1=\sum_{n\le x}\sum_{\alpha|n\atop
\alpha|d}\mu(\alpha)=\sum_{\alpha|d}\mu(\alpha)[{x\over \alpha}]
={\varphi(d)\over d}x+O(2^{\omega(d)}).
\end{equation}
By the principle of inclusion and exclusion we find that
$$S_d(x)=\sum_{m\le \sqrt{x}\atop (d,m)=1}\mu(m)T_d({x\over m^2}).$$
Hence, on invoking (\ref{teedee}), we find
$$S_d(x)=x{\varphi(d)\over d}\sum_{m\le \sqrt{x}\atop (d,m)=1}
{\mu(m)\over m^2}+O(2^{\omega(d)}\sqrt{x}).$$
and hence, on completing the sum,
$$S_d(x)=x{\varphi(d)\over d}\sum_{m=1\atop (d,m)=1}^{\infty}{\mu(m)\over m^2}+
O(2^{\omega(d)}\sqrt{x})$$
Note that
$$\sum_{m=1\atop (d,m)=1}^{\infty}{\mu(m)\over m^2}=\prod_{p\nmid d}(1-{1\over p^2})
={1\over \zeta(2)\prod_{p|d}(1-1/p^2)}.$$
Using this and (\ref{startphi}) the proof is completed. \qed\\

Let $d(n)$ denote the number of divisors of $n$. 
We have $2^{\omega(n)}\le d(n)$ with equality iff $n$ is squarefree. The
estimates below also hold with $2^{\omega(n)}$ replaced by $d(n)$.
\begin{lem}
We have 
$$\sum_{d\le x}^{\infty}{2^{\omega(d)}\over d^{3/2}}=O(1),~ 
\sum_{d\le x}{2^{\omega(d)}\over \sqrt{d}}=O(\sqrt{x}\log x),~
\sum_{d\le x}{4^{\omega(d)}\over d}=O(\log ^3 x).$$
\end{lem}
{\it Proof}. 
Using the convergence of $\sum_p p^{-3/2}$ we find by (\ref{pprod}) that
$\sum_{d=1}^{\infty}2^{\omega(d)}d^{-3/2}=O(1)$.
The remaining estimates follow on invoking Theorem 1 at p. 201 of 
Tenenbaum's book \cite{tenenbaum} together with partial integration. \qed

\subsection{Proof of Theorem \ref{main}}
Note that
$$C_1(x)=\sum_{a\le x}\sum_{b\le x}\mu(a)^2\sum_{d|a,~d|b}\mu(d)=
\sum_{d\le x}\mu(d)\sum_{a\le x\atop d|a}\mu(a)^2\sum_{b_1\le x/d}1,$$
after swapping the summation order. Using $[x/d]=x/d+O(1)$, we then obtain
$$C_1(x)=x\sum_{d\le x}{\mu(d)\over d}\sum_{a\le x\atop d|a}\mu(a)^2
+O(x\log x).$$
On noting that
\begin{equation}
\label{simpel}
\sum_{a\le x\atop d|a}\mu(a)^2=\mu(d)^2\sum_{n\le x/d\atop (d,n)=1}\mu(n)^2
=\mu(d)^2S_d({x\over d})
\end{equation}
and $\mu(d)=\mu(d)^3$, we find
$$C_1(x)=x\sum_{d\le x}{\mu(d)\over d}S_d({x\over d})+O(x\log x).$$
On using  Lemma 1 we obtain the estimate
$$C_1(x)={x^2\over \zeta(2)}\sum_{d\le x}{\mu(d)\over d^2\prod_{p|d}(1+1/p)}
+O(\sqrt{x}\sum_{d\le x}{2^{\omega(d)}\over \sqrt{d}})+O(x\log x).$$
On completing the latter sum and noting that
$$\sum_{d=1}^{\infty}{\mu(d)\over d^2\prod_{p|d}(1+1/p)}
=\prod_p\Big(1-{1\over p(p+1)}\Big),$$
we obtain
$$C_1(x)={x^2\over \zeta(2)}\prod_p\Big(1-{1\over p(p+1)}\Big)
+O(\sqrt{x}\sum_{d\le x}{2^{\omega(d)}\over \sqrt{d}})+O(x\log x).$$
Estimate (1) now follows on invoking Lemma 2.\\
\indent The proof of (\ref{twee}) is very
similar to the proof of (\ref{een}). We start by noting that
$$C_2(x)=\sum_{a\le x}\sum_{b\le x}\mu(a)^2\mu(b)^2\sum_{d|a,~d|b}\mu(d).$$
On swapping the summation order, we obtain
\begin{equation}
\label{inex-1}
C_2(x)=\sum_{d\le x}\mu(d)\sum_{a\le x\atop d|a}\mu(a)^2\sum_{b\le x\atop
d|b}\mu(b)^2.
\end{equation}
On noting that $\mu(d)=\mu(d)^5$ and invoking (\ref{simpel}) we obtain
\begin{equation}
\label{inex}
C_2(x)=\sum_{d\le x}\mu(d)S_d({x\over d})^2.
\end{equation}
On using Lemma 1 we obtain the estimate
$$C_2(x)={x^2\over \zeta(2)^2}\sum_{d\le x}{\mu(d)\over d^2\prod_{p|d}(1+1/p)^2}
+O(x^{3/2}\sum_{d\le x}{2^{\omega(d)}\over d^{3/2}})
+O(x\sum_{d\le x}{4^{\omega(d)}\over d}).$$
On completing the first sum and noting that
$$\sum_{d=1}^{\infty}{\mu(d)\over d^2\prod_{p|d}(1+1/p)^2}
=\prod_p\Big(1-{1\over (p+1)^2}\Big),$$
we find
$$C_2(x)={x^2\over \zeta(2)^2}\prod_p\Big(1-{1\over (p+1)^2}\Big)
+O(x^{3/2}\sum_{d\le x}{2^{\omega(d)}\over d^{3/2}})
+O(x\sum_{d\le x}{4^{\omega(d)}\over d}).$$
On invoking Lemma 2 
estimate (\ref{twee}) is then established. \qed
\subsection{Proof of Theorem \ref{main2}}
\label{st2}
We write $[n,m]$ for the least common multiple of $n$ and $m$,
and $(n,m)$ for the greatest common divisor. Recall that $(n,m)[n,m]=nm$.\\
\indent Note that
$$I_3(x)=\sum_{a,b,c\le x}\sum_{d_1|a\atop d_1|b}\mu(d_1)
\sum_{d_2|a\atop d_2|c}\mu(d_2)\sum_{d_3|b\atop d_3|c}
\mu(d_3),$$
which can be rewritten as
$$I_3(x)=\sum_{{[d_1,d_2]\le x\atop [d_1,d_3]\le x}\atop [d_2,d_3]\le
x}\mu(d_1)\mu(d_2)\mu(d_3)[{x\over [d_1,d_2]}]
[{x\over [d_1,d_3]}][{x\over [d_2,d_3]}].$$
Now put
$$J_1(x)=\sum_{{[d_1,d_2]\le x\atop [d_1,d_3]\le x}
\atop [d_2,d_3]\le x}{\mu(d_1)\mu(d_2)\mu(d_3)\over
[d_1,d_2][d_1,d_3][d_2,d_3]},~J_2(x)=\sum_{{[d_1,d_2]\le x\atop [d_1,d_3]\le x}
\atop [d_2,d_3]\le x}{1\over [d_1,d_2][d_1,d_3]},$$ 
$$J_3(x)=\sum_{{[d_1,d_2]\le x\atop [d_1,d_3]\le x}
\atop [d_2,d_3]\le x}{1\over [d_1,d_2]}{\rm ~and~}J_4(x)=\sum_{{[d_1,d_2]\le x\atop [d_1,d_3]\le x}
\atop [d_2,d_3]\le x}1.$$
Using that $[x]=x+O(1)$ we find that
\begin{equation}
\label{i3}
I_3(x)=x^3J_1(x)+O(x^2J_2(x))+O(xJ_3(x))+O(J_4(x)).
\end{equation}
We will show first that
$$J_1(x)=\sum_{d_1=1}^{\infty}\sum_{d_2=1}^{\infty}
\sum_{d_3=1}^{\infty}{\mu(d_1)\mu(d_2)\mu(d_3)\over
[d_1,d_2][d_1,d_3][d_2,d_3]}+O\Big({\log x\over x}\Big).$$
To this end it is enough, by symmetry of the argument of the
sum, to show that
\begin{equation}
\label{symorder}
\sum_{[d_1,d_2]> x}\sum_{d_3\ge 1}{1\over
[d_1,d_2][d_1,d_3][d_2,d_3]}=O\Big({\log x\over x}\Big).
\end{equation}
Put $(d_1,d_2)=\alpha$, $(d_1,d_3)=\beta$ and $(d_2,d_3)=
\gamma$. 
Since $\alpha|d_1$ and $\beta|d_1$, we can write
$d_1=[\alpha,\beta]\delta_1$ for some integer $\delta_1\ge 1$,
and similarly $d_2=[\alpha,\gamma]\delta_2$, $d_3=[\beta,\gamma]\delta_3$. Note
that any triple 
$(d_1,d_2,d_3)$ corresponds to a uniqe 6-tuple $(\alpha,\beta,\gamma,\delta_1,\delta_2,\delta_3)$.
Since $\alpha(\delta_1,\delta_2)$ divides $([\alpha,\beta]\delta_1,[\alpha,\gamma]\delta_2)$
on the one hand and $([\alpha,\beta]\delta_1,[\alpha,\gamma]\delta_2)=(d_1,d_2)=\alpha$
on the other, it follows that
$(\delta_1,\delta_2)=1$ and likewise $(\delta_1,\delta_3)=(\delta_2,\delta_3)=1$. 
Write $u=\alpha\beta\gamma/(\alpha,\beta,\gamma)^2$.
On noting that $((d_1,d_2),(d_2,d_3))=(d_1,d_2,d_3)=((d_1,d_2),(d_1,d_3),(d_2,d_3))$ we infer
that $(\alpha,\beta)=(\alpha,\gamma)=(\beta,\gamma)=(\alpha,\beta,\gamma)$
and hence we find that $[d_1,d_2]=u\delta_1\delta_2$, $[d_1,d_3]=u\delta_1\delta_3$
and $[d_2,d_3]=u\delta_2\delta_3$. Now
$$\sum_{[d_1,d_2]> x}\sum_{d_3\ge 1}{1\over
[d_1,d_2][d_1,d_3][d_2,d_3]}\le \sum_{\alpha,\beta,\gamma}{1\over u^3}
\sum_{\delta_1\delta_2>x/u}\sum_{\delta_3\ge 1}{1\over (\delta_1\delta_2
\delta_3)^2},$$
where the triple sum is over all 6-tuples $(\alpha,\beta,\gamma,\delta_1,\delta_2,\delta_3)$ and is of order
$$O\Big(\sum_{\alpha,\beta,\gamma}{1\over u^3}
\sum_{\delta_1 \delta_2> x/u}{1\over (\delta_1\delta_2)^2}
\Big)
=O\Big(\sum_{\alpha,\beta,\gamma}{1\over u^3}\sum_{n>x/u}{d(n)\over n^2}\Big)
=O\Big({\log x\over x}\sum_{\alpha,\beta,\gamma}{1\over u^2}
\Big),$$
where we used the well-known estimate $\sum_{n>x}d(n)n^{-2}=O(\log x /x)$.
Now
\begin{equation}
\label{lunch}
\sum_{\alpha,\beta,\gamma}{1\over u^2}
=\sum_{\alpha,\beta,\gamma}
{(\alpha,\beta,\gamma)^4\over (\alpha \beta \gamma)^2}
=O\Big(\sum_{d=1}^{\infty}
{1\over d^2}\sum_{\alpha ',\beta ',\gamma '}{1\over
(\alpha '\beta '\gamma ')^2}\Big)=O(1),
\end{equation}
where we have written $(\alpha,\beta,\gamma)=d$, $\alpha=d\alpha'$, 
$\beta=d\beta'$ and $\gamma=d\gamma'$.
Thus we have established equation (\ref{symorder}).\\
\indent In the same vein $J_2(x)$ can be estimated to be
\begin{eqnarray}
J_2(x)&=&O\Big(\sum_{\alpha,\beta,\gamma}\sum_{\delta_1\delta_2\le x/u\atop 
{\delta_1\delta_3\le x/u\atop \delta_2\delta_3\le x/u}}{1\over [d_1,d_2][d_1,d_3]}\Big)
=O\Big(\sum_{\alpha,\beta,\gamma}{1\over u^2}
\sum_{\delta_1\delta_2\delta_3\le (x/u)^{3/2}}{1\over \delta_1^2
\delta_2\delta_3}\Big)\cr
&=&O\Big(\sum_{\alpha,\beta,\gamma}{1\over u^2}\sum_{\delta_2\delta_3\le (x/u)^{3/2}}
{1\over \delta_2\delta_3}\Big)=
O\Big(\sum_{\alpha,\beta,\gamma}{1\over u^2}\sum_{n\le (x/u)^{3/2}}
{d(n)\over n}\Big)\nonumber
\end{eqnarray}
Using the classical estimate $\sum_{n\le x}d(n)/n=O(\log ^2 x)$ and (\ref{lunch}), one obtains
$J_2(x)=O(\log^2 x)$.\\ 
\indent Note that $0\le J_4(x)\le xJ_3(x)\le x^2J_2(x)$.
Using (\ref{i3}) we see that it remains to evaluate
the triple infinite sum, which we rewrite as
$$\sum_{d_1=1}^{\infty} \sum_{d_2=1}^{\infty}
\sum_{d_3=1}^{\infty}{\mu(d_1)\mu(d_2)\mu(d_3)(d_1,d_2)(d_1,d_3)
(d_2,d_3)\over (d_1d_2d_3)^2},$$
which can be rewritten as
$$\sum_{d_1=1}^{\infty}{\mu(d_1)\over d_1^2}
\sum_{d_2=1}^{\infty}{\mu(d_2)(d_1,d_2)\over d_2^2}
\sum_{d_3=1}^{\infty}{\mu(d_3)(d_1,d_3)(d_2,d_3)\over d_3^2}.$$
Note that the argument of the inner sum is multiplicative
in $d_3$. By Euler's product identity (\ref{pprod}) it is zero if $(d_1,d_2)>1$
and $\zeta(2)^{-1}\prod_{p|d_1d_2}(1+1/p)^{-1}$ otherwise.
Thus the latter triple sum is seen to yield
$${1\over \zeta(2)}\sum_{d_1=1}^{\infty}{\mu(d_1)\over
d_1^2\prod_{p|d_1}(1+1/p)}\sum_{d_2=1\atop (d_1,d_2)=1}
^{\infty}{\mu(d_2)\over d_2^2\prod_{p|d_2}(1+1/p)},$$
the argument of the inner sum is multiplicative in $d_2$
and proceeding as before we 
obtain that it equals
$${1\over \zeta(2)}\prod_p\Big(1-{1\over p(p+1)}\Big)\sum_{d_1=1}^{\infty}
{\mu(d_1)\over d_1^2\prod_{p|d_1}(1+{1\over p})\prod_{p|d_1}(1-{1\over p(p+1)})},$$
which is seen to equal
$${1\over \zeta(2)}\prod_p\Big(1-{2\over p(p+1)}\Big),$$
which by equation (\ref{trafo}) equals $K_2$. \qed

\section{Numerical aspects}
Direct evaluation of the constants $K_1$ and $K_2$ 
through (3), respectively (4) yields only about five decimal digits of 
precision. 
By
expressing $K_1$ and $K_2$ as infinite products involving 
$\zeta(k)$ for $k\ge 2$, they can be computed with high
precision. To this end Theorem 1 of \cite{Mor} can be used. 
The error analysis can be dealt with using Theorem 2 of \cite{Mor}. Using 
\cite[Theorem 1]{Mor}  it is inferred
that $$K_1=\prod_{k\ge 2}\zeta(k)^{-e_k},{\rm ~where~}e_k
={\sum_{d|k}b_d\mu({k\over d})\over k}\in \Bbb Z,$$
with the sequence $\{b_k\}_{k=0}^{\infty}$ defined by $b_0=2$
and $b_1=-1$ and $b_{k+2}=-b_{k+1}+b_k$. Using
the same theorem, it is seen that 
$$K_2={1\over 2}\prod_{k\ge 2}\{\zeta(k)(1-2^{-k})\}^{-f_k}, 
{\rm ~where~}f_k={\sum_{d|k}(-2)^d\mu({k\over d})\over k}\in \Bbb Z.$$
\indent Typically in analytic number theory constants
of the form $\prod_p f(1/p)$ with $f$ rational arise as densities. 
Their numerical evaluation was considered by the author in \cite{Mor}.
By similar methods any constant of the form $\prod_p f(1/p)$ with $f$ an 
analytic function on the unit disc satisfying $f(0)=1$ and $f'(0)=0$ can 
be evaluated \cite{MP}.

\section{Related problems}
Let us call a couple $(a,b)$ with $a,~b\le x$, $a$ and $b$
coprime and either $a$
or $b$ squarefree, {\it weakly carefree}. A little thought reveals
that
$C_3(x)=2C_1(x)-C_2(x)$. By Theorem 1 it then follows that the probability
$K_3$ that a couple is weakly carefree equals
$K_3=2K_1-K_2\approx 0.5697515829$.\\
\indent The problem of
estimating $I_3(x)$
has the following natural generalisation.
Let $k\ge 2$ be an integer and let $I_k(x)$ be the number of $k$-tuples
$(a_1,\ldots,a_k)$ with $1\le a_i\le x$ for $1\le i\le k$ such
that $(a_i,a_j)=1$ for every $1\le i\ne j\le k$. The number of
$k$-tuples such that none of the gcd's is divisible by some
fixed prime $p$ is easily seen to be
$$\sim x^k\Big(\Big(1-{1\over p}\Big)^k+{k\over p}\Big(1-{1\over p}\Big)^{k-1}\Big)=
x^k\Big(1-{1\over p}\Big)^{k-1}\Big(1+{k-1\over p}\Big).$$
Thus, it seems plausible that
\begin{equation}
\label{vermoeden}
I_k(x)\sim x^k\prod_p \Big(1-{1\over p}\Big)^{k-1}\Big(1+{k-1\over p}\Big),~~~~~(x\rightarrow \infty).
\end{equation}
For $k=2$ and $k=3$ (by Theorem 2 and equation (\ref{trafo}) this is true. 
In 2000 I did not see how to prove this for arbitrary $k$, however the conjecture
(\ref{vermoeden}) was established soon afterwards (in 2002) by 
L. T\'oth \cite{toth}, who proved that for $k\ge 2$ we have
\begin{equation}
\label{tottie}
I_k(x)=x^k\prod_{p}\Big(1-{1\over p}\Big)^{k-1}\Big(1+{k-1\over p}\Big)+O(x^{k-1}\log^{k-1}x).
\end{equation}
\indent Let $I_k^{(u)}(x)$ denote the number of $k$-tuples $(a_1,\ldots,a_k)$ with
$1\le a_i\le x$ that are pairwise coprime and moreover satisfy $(a_i,u)=1$ for
$1\le i\le k$. It is easy to see that
$$I_{k+1}^{(u)}(n)=\sum_{j=1\atop (j,u)=1}^n I_k^{(ju)}(n).$$
Note that $I_1^{(u)}(n)=T_u(n)$ can be estimated by (\ref{teedee}). Then by recursion
with respect to $k$ an estimate for $I_k^{(u)}(n)$ can be 
established that implies (\ref{tottie}).\\
\indent In \cite{HM} Havas and Majewski considered the problem of counting
the number of $n$-tuples of natural numbers that are pairwise not coprime. They
suggested that the density $\delta_n$ of these tuples should be
\begin{equation}
\label{simpel2}
\delta_n=\Big(1-{1\over \zeta(2)}\Big)^{({n\atop 2})}.
\end{equation}
The probability that a pair of integers is not coprime is $1-1/\zeta(2)$. Since
there are $({n\atop 2})$ pairs of integers in an $n$-tuple, one might naively
espect the probability for this problem to be as given by (\ref{simpel2}). \\
\indent T. Freiberg \cite{F} studied this problem for $n=3$ using
my approach to estimate $I_3(x)$ 
(it seems that the recursion method of T\'oth cannot be applied here). 
Freiberg showed that the density of triples 
$(a,b,c)$ with $(a,b)>1,~(a,c)>1$ and $(b,c)>1$ equals
$$F_3=1-{3\over \zeta(2)}+3K_1-K_2\approx 0.1742197830347247005,$$
whereas $(1-1/\zeta(2))^3\approx 0.06$. 
Thus the guess of Havas and Majewski for $n=3$ is false. Indeed, it is easy to see 
(as Peter Pleasants pointed out to the author \cite{pleasants}) 
that for every $n\ge 3$ their guess is false. Since all $n$-tuples of even numbers
are pairwise not coprime, $\delta_n$, if it exists, satisfies $\delta_n\ge 2^{-n}$.
Since $({n\atop 2})\ge n$ and $1-1/\zeta(2)<0.4$ the predicted density by
Havas and Majewski \cite{HM} satisfies $\delta_n<2^{-n}$ for $n\ge 3$ and so must be false.\\
\indent In 2006 the author learned \cite{keith} that the result of Freiberg
is implicit in the PhD thesis of R.N. Buttsworth \cite{bob} and indeed can be found there in
more general form. Buttsworth showed that the density of 
relatively prime $m$-tuples for which $k$ prescribed $(m-1)$-tuples have gcd 1
equals $Z_k(m)$ given in (\ref{zkm}).
Consequently by inclusion and exclusion the set of relatively prime $m$-tuples such that every
$(m-1)$-tuple fails to be
relatively prime  has density
$$\sum_{k=0}^m  (-1)^k {m\choose k} Z_k(m).$$
For $m=3$ this yields  $1/\zeta(3) -3/\zeta(2) +3K_1 - K_2$.
So the density of relatively prime 3-tuples such that at  least
one 2-tuple is relatively prime, is equal to $3/\zeta(2) - 3K_1 + K_2$.
However this is also equal to the  density of 3-tuples such that at  least
one 2-tuple is relatively prime.
Hence the  density of 3-tuples such that  all 2-tuples are not
relatively prime is
$1 -3/\zeta(2) +3K_1 - K_2$, 
which is Freiberg's formula.\\
\indent To close this discussion, we like to remark that Freiberg established his result with
error term $O(x^2 \log^2 x)$ and that Buttsworth's result gives only a density.\\
\indent Some related open problems are as follows: 
\begin{Problem}$~~$\\
{\rm a)} To compute the density of $n$-tuples such that at least $k$ pairs are coprime.\\
{\rm b)} To compute the density of $n$-tuples such that exactly $k$ pairs are coprime.
\end{Problem}
\begin{Problem}$~~$\\
To compute the density of $n$-tuples such that all pairs are not coprime.
\end{Problem}
{\tt Remark}. Recently Jerry Hu \cite{Hu} announced that he solved Problem 1.

\section{Conclusion}
In stark constrast to what experience from daily life suggests,
(strongly) carefree couples are quite common...\\

\noindent {\tt Acknowledgement}. The author likes to thank Steven Finch for 
bringing Schroeder's conjecture to his attention and his
instignation to write down these results. 
Also Finch and de Weger pointed out that
one has $\sum_{n\le x}k(n)=\zeta(2)K_1x^2/2+O(x^{3/2})$, where
$k(n)=\prod_{p|n}p$, and that in \cite{Morun} 
the $K_1$ was inadvertently dropped. For a proof of this formula see Eckford Cohen 
\cite[Theorem 5.2]{Cohen}.\\
\indent The author likes to thank Tristan Freiberg 
and Jerry Hu for pointing out some references and helpful comments. Keith Matthews
provided me kindly with very helpful information concerning the relevant results of his
former PhD student Buttsworth. In particular he pointed out how Freiberg's result follows
from that of Buttsworth.

\medskip\noindent {\footnotesize Max-Planck-Institut f\"ur Mathematik,
Vivatsgasse 7, D-53111 Bonn, Germany.\\
e-mail: {\tt moree@mpim-bonn.mpg.de}\\
\end{document}

\end{document}

\begin{eqnarray}
\label{teedee}
T_d(x)&=&\sum_{n\le x\atop (n,d)=1}1=\sum_{n\le x}\sum_{\alpha|n\atop
\alpha|d}\mu(\alpha);\cr
&=&\sum_{\alpha|d}\mu(\alpha)[{x\over \alpha}]={\varphi(d)\over d}x+O(2^{\omega(d)}).
\end{eqnarray}